\documentclass[reqno,12pt,draft]{amsart}
\usepackage{amsmath,amsthm,amsfonts,amssymb}
\date{}
\newtheorem{theorem}{Theorem}
\newtheorem{lemma}{Lemma}

\begin{document}

\centerline{\large{\sc Gradient bounds for solutions of
elliptic}}

\centerline{\large{\sc
 and parabolic equations}}

\vskip .1in

\centerline{{\sc Vladimir~I.~Bogachev$^{a}$, Giuseppe~Da~Prato$^{b}$,
Michael~R\"ockner$^{c}$,}}

\centerline{{\sc and Zeev~Sobol$^{d}$}}

\vskip .1in

$^{a}$:
{\footnotesize Department of Mechanics and Mathematics,
         Moscow State University, 119992 Moscow, Russia}

$^{b}$:
{\footnotesize
Scuola Normale Superiore di Pisa, Piazza dei Cavalieri 7,
I--56125 Pisa, Italy}

$^{c}$:
{\footnotesize
Fakult\"at f\"ur Mathematik,
Universit\"at Bielefeld, D--33501 Bielefeld, Germany}

$^{d}$:
{\footnotesize
Department of Mathematics, University of Wales Swansea,
Singleton Park, Swansea SA2 8PP, UK}
\vskip .2in

\centerline{{\bf Abstract}}
\vskip .1in

{\small
Let $L$ be a second order elliptic operator on $\mathbb{R}^d$ with
a constant diffusion matrix and a dissipative (in a weak sense) drift
$b\in L^p_{loc}$ with some $p>d$.
We assume that $L$ possesses a Lyapunov function,
but no local boundedness of $b$ is assumed. It is known that
then there exists a unique probability measure $\mu$ satisfying
the equation $L^*\mu=0$ and that the closure of $L$ in
$L^1(\mu)$ generates a
Markov semigroup $\{T_t\}_{t\ge 0}$ with the
resolvent $\{G_\lambda\}_{\lambda > 0}$. We prove
that, for any Lipschitzian function
$f\in L^1(\mu)$ and all $t,\lambda>0$,
the functions $T_tf$ and $G_\lambda f$ are Lipschitzian and
$$
|\nabla T_tf(x)|\leq T_t|\nabla f|(x)
\quad\hbox{and}\quad
|\nabla G_\lambda f(x)|\leq
\frac{1}{\lambda} G_\lambda |\nabla f|(x).
$$
An analogous result is proved in the
parabolic case.
}

\vskip .3in

Suppose that for every $t\in [0,1]$, we are given a
a strictly positive definite symmetric
matrix
$A(t)=(a^{ij}(t))$
and
a measurable vector field
$x\mapsto b(t,x)=(b^1(t,x),\ldots,b^n(t,x))$.

Let $L_t$ be the elliptic operator on $\mathbb{R}^d$ given by
\begin{equation}\label{op}
L_tu(x)=\sum\limits_{i,j\le d}
a^{ij}(t,x)\partial_{x_i}\partial_{x_j}u(x)
+\sum\limits_{i\le d} b^i(t,x)\partial_{x_i} u(x).
\end{equation}

Suppose that $A$ and $b$ satisfy the following hypotheses:

\begin{enumerate}
\item[(Ha)]
$\sup_{t\in [0,1]}\bigl( \|A(t)\|+\|A(t)^{-1}\|\bigr)<\infty$,
$\sup_{t\in [0,1]}
\|b(t,\,\cdot\,)\|_{L^p(U)}<\infty$ for every ball
$U$ in $\mathbb{R}^d$ with some $p>d$, $p\geq2$.

\item[(Hb)] $b$ is {\it dissipative} in the following sense:
for every $t\in [0,1]$ and every
$h\in \mathbb{R}^d$, there exists a measure zero set
$N_{t,h}\subset\mathbb{R}^d$ such that
$$
\bigl(b(t,x+h)-b(t,x),h\bigr)\le 0\quad
\hbox{ for all
$x\in\mathbb{R}^d\setminus N_{t,h}$.}
$$

\item[(Hc)]
for every $t\in [0,1]$,
there exists a {\it Lyapunov function} $V_t$ for $L_t$, i.e., a
nonnegative $C^2$-function $V_t$ such that
$V_t(x)\to+\infty$ and $L_tV_t(x)\to -\infty$ as $|x|\to\infty$.
\end{enumerate}

We consider the parabolic equation
\begin{equation}\label{par1}
\frac{\partial u}{\partial t}=L_tu,\quad u(0,x)=f(x),
\end{equation}
where $f$ is a bounded Lipschitz function. A locally integrable
function $u$
on $[0,1]\times\mathbb{R}^d$ is called a solution if,
for every $t\in (0,1]$, one has $u(t,\,\cdot\,)\in
W^{1,2}_{loc}(\mathbb{R}^d)$, the functions
$\partial_{x_i}\partial_{x_j}u$ and
$b^i\partial_{x_i}u$ are integrable on the sets
$[0,1]\times K$ for every cube $K$ in
$\mathbb{R}^d$,
and for every $\varphi\in C_0^\infty(\mathbb{R}^d)$ and all
$t\in [0,1]$ one has
$$
\int_{\mathbb{R}^d} u(t,x)\varphi(x)\, dx
=\int_{\mathbb{R}^d} f(x)\varphi(x)\, dx
+\int_0^t\int_{\mathbb{R}^d}
L_s\varphi(x)\, u(s,x)\, dx ds.
$$

In the case where $A$ and $b$ are independent of~$t$,
     so that we have a single operator~$L$,
Hypotheses (Ha) and (Hc) imply
(see \cite{BR} and \cite{BRS})
that there exists a unique probability
measure $\mu$ on $\mathbb{R}^d$
such that $\mu$ has a
strictly positive continuous weakly differentiable
density $\varrho$, $|\nabla\varrho|\in L^p_{loc}(\mathbb{R}^d)$,
and $L^*\mu=0$ in the following
weak sense:
$$
\int Lu\,d\mu=0\quad \hbox{for all $u\in C_0^\infty(\mathbb{R}^d)$.}
     $$
The closure
$\overline{L}$ of $L$ with domain $C_0^\infty(\mathbb{R}^d)$
in $L^1(\mu)$
generates a Markov
semigroup $\{T_t\}_{t\ge 0}$ for which $\mu$ is invariant.
Let
$D(\overline{L})$ denote the domain of $\overline{L}$ in
     $L^1(\mu)$  and let
$\{G_\lambda\}_{\lambda > 0}$ denote the corresponding
resolvent, i.e., $G_\lambda=(\lambda -\overline{L})^{-1}$.
The restrictions of $T_t$ and $G_\lambda$
to $L^2(\mu)$ are contractions on~$L^2(\mu)$.
In particular, if $v\in D(\overline{L})$ is such that
$\lambda v-\overline{L}v=g\in L^2(\mu)$, then
$v\in L^2(\mu)$. Moreover, it follows by \cite[Theorem~2.8]{BRS}
that one has
     $v\in H_{loc}^{2,2}(\mathbb{R}^d)$ and
     $\overline{L}v=Lv$~a.e., so that one has~a.e.
     \begin{equation}\label{e1}
     \lambda v-Lv=g.
\end{equation}
In fact, due to our assumptions
on the coefficients of $L$ one has even
$v\in W^{p,2}_{loc}(\mathbb{R}^d)$ (see \cite{Chicco}).
It has been shown in \cite{BDRS} that for every
     function $f\in L^1(\mu)$ that is Lipschitzian with
constant $C$ and all $t,\lambda>0$, the continuous
     version of the function $T_tf$ is Lipschitzian with constant~$C$,
     and the continuous version of $G_\lambda f$ is Lipschitzian
     with constant~$\lambda^{-1}C$. Here we establish pointwise
     estimates in both cases and prove their parabolic analogue.
The main results of this work are the following two theorems.

\begin{theorem}
Suppose that $A$ and $b$ are independent of $t$ and
satisfy {\rm(Ha)}, {\rm(Hb)} and
{\rm(Hc)}. Then, for any Lipschitzian function
$f\in L^1(\mu)$ and all $t,\lambda>0$,
$T_tf$ and $G_\lambda f$ have Lipschitzian versions such that
\begin{equation}\label{local}
|\nabla T_tf(x)|\leq T_t|\nabla f|(x)
\quad\hbox{and}\quad
|\nabla G_\lambda f(x)|\leq \frac{1}{\lambda}
G_\lambda |\nabla f|(x)
\end{equation}
for the corresponding continuous versions.
In particular,
\begin{equation}\label{glob}
\sup\limits_{x,t}|\nabla T_tf(x)|\leq \sup\limits_x|\nabla f(x)|,
\quad
\sup\limits_{x}|\nabla G_\lambda f(x)|\leq \frac{1}{\lambda}
\sup\limits_x|\nabla f(x)|.
\end{equation}
\end{theorem}

\begin{theorem}
Suppose that $A$ and $b$
satisfy {\rm(Ha)}, {\rm(Hb)} and
{\rm(Hc)}. Then, for any bounded Lipschitzian function
$f$ there is a solution $u$ of equation {\rm(\ref{par1})} such
that for all $t$ one has
\begin{equation}\label{glob2}
\sup\limits_{x}|\nabla u(t,x)|\leq \sup\limits_x|\nabla f(x)| .
\end{equation}
\end{theorem}

In the case where $A=I$ and $b=0$, estimate (\ref{glob2})
     has been established in \cite{Kahane1}, \cite{Kahane2}
for solutions of boundary problems in bounded domains.
It should be noted that gradient estimates of the type
     $$
\sup\limits_x |\nabla u(x,t)|\le C(t)\sup\limits_x |f(x)|
$$
for solutions of parabolic equations have been obtained
     by many authors, see, e.g.,
\cite{BF}, \cite{BFL},
\cite{FMP}, \cite{Wang}, and the references therein. Such estimates do not
     require (Hb) and  one  has  $C(t)\to  +\infty$  as  $t\to  0$  or
$t\to +\infty$. In contrast to this type of estimates, our theorems
mean a contraction property on Lipschitz functions rather than a
smoothing property.
     It is likely that some  results of the cited works,
     established for sufficiently regular~$b$,
can be extended to more general drifts satisfying just~(Ha),
     but not~(Hb).

A short proof of the following result can be found in
\cite{BDRS}.

\begin{lemma}
Suppose that $b$ is infinitely differentiable, Lipschitzian,
and strongly dissipative, so for some
$\alpha>0$, one has
$$
\bigl(b(x+h)-b(x),h\bigr)\le -\alpha(h,h)
\quad \hbox{for all $x,h\in\mathbb{R}^d$.}
$$
Then, for any $\lambda>0$ and any smooth bounded Lipschitzian
function~$f$, one has pointwise
$$
|\nabla G_\lambda f|\le G_\lambda |\nabla f| .
$$
In particular,
$\sup\limits_x|\nabla G_\lambda f(x)|\leq
\lambda^{-1}\sup\limits_x|\nabla f(x)|$.
\end{lemma}

\begin{proof}[Proof of Theorem~1]
The estimate with the suprema has been proven
     in \cite{BDRS}, and the stronger pointwise estimate can be
     derived from that proof. For the reader's convenience,
     instead of recursions to the steps of the proof in \cite{BDRS}
     we reproduce the whole proof and explain why it yields a
     stronger conclusion.
We recall that if a sequence of functions on $\mathbb{R}^d$ is
uniformly Lipschitzian with constant $L$ and
bounded at a point, then it contains a
subsequence that converges uniformly on every ball to a function
that is Lipschitzian with the same constant. Therefore,
approximating $f$ in $L^1(\mu)$ by a sequence of bounded smooth
functions $f_j$ with
$$
\sup_x |\nabla f_j(x)|\le \sup_x |\nabla f(x)|,
$$
it suffices to prove (\ref{glob}) for smooth bounded~$f$.
Moreover, due to Euler's formula
$T_tf=\lim_n\Big(\frac{t}{n}G_{\frac{t}{n}}\Big)^nf$, it suffices
to establish the resolvent estimate.
First we construct a suitable sequence of smooth strongly dissipative
Lipschitzian vector fields $b_k$  such that $b_k\to b$ in
$L^p(U,\mathbb{R}^d)$ for every ball $U$ as $k\to\infty$.
Let $\sigma_j(x)=j^{-d}\sigma(x/j)$, where $\sigma$ is
a smooth compactly supported probability density. Let
$\beta_{j}:=b *\sigma_j$.
Then $\beta_j$ is smooth and dissipative and $\beta_j\to b$,
$j\to\infty$, in
$L^p(U,\mathbb{R}^d)$ for every ball~$U$.
For every $\alpha>0$,
the mapping $I-\alpha \beta_j$ is a
homeomorphism of $\mathbb{R}^d$ and the inverse mapping
$(I-\alpha \beta_j)^{-1}$
is Lipschitzian with constant~$\alpha^{-1}$ (see \cite{Brezis}).
Let us consider the Yosida approximations
$$
F_\alpha(\beta_j):=\alpha^{-1}\bigl((I-\alpha \beta_j)^{-1}-I\bigr)
=\beta_j\circ(I-\alpha \beta_j)^{-1}.
$$
It is known (see \cite[Ch.~II]{Brezis})
that $|F_\alpha(\beta_j)(x)|\le |\beta_j(x)|$,
the mappings $F_\alpha(\beta_j)$
converge locally uniformly to $\beta_j$
as $\alpha\to 0$, and one has
$$
\bigl(F_\alpha(\beta_j)(x)-F_\alpha(\beta_j)(y),x-y\bigr)\le 0.
$$
Thus, the
sequence $b_k:=F_{\frac{1}{k}}(b*\sigma_k) -\frac{1}{k} I$,
$k\in\mathbb{N}$, is the desired one.
For every $k\in\mathbb{N}$, let $L_k$ be the elliptic operator defined
by~(\ref{op}) with the same constant matrix $A$
and drift~$b_k$ in place of~$b$.
Let
$\mu_k=\varrho_k\, dx$ be the corresponding invariant probability
measure and let
$G^{(k)}_\lambda$ denote the associated resolvent family
on~$L^1(\mu_k)$.
Since $b_k$ is smooth, Lipschitzian and strongly dissipative,
$v_k:=G^{(k)}_\lambda f$ is smooth, bounded, Lipschitzian and
$$
\sup\limits_x|v_k(x)|\leq \frac{1}{\lambda}\sup\limits_x|f(x)|
\quad\hbox{and}\quad
\sup\limits_x|\nabla v_k(x)|\leq
\frac{1}{\lambda}\sup\limits_x|\nabla f(x)|
$$
by the lemma. Moreover, for every ball~$U\subset \mathbb{R}^d$,
the functions $v_k$ are uniformly bounded in the Sobolev space
$W^{2,2}(U)$, since the mappings $|b_k|$ are bounded in $L^p(U)$
uniformly in $k$ and $f$ is bounded.
This follows from the fact that for any solution
$w\in W^{2,2}(U)$ of the equation
$\sum\limits_{i,j\le d} a^{ij}\partial_{x_i}\partial_{x_i}w+
\sum\limits_{i\le d} b^{i}\partial_{x_i}\partial_{x_i}w-\lambda w=g$
one has $\|w\|_{W^{2,2}(U)}\le C\|w\|_{L^2(U)}$, where
$C$ is a constant that depends on $U$, $A$, and the quantity
$\kappa:=\|g\|_{L^2(U)}+\| |b|\|_{L^p(U)}$ in such a way that
as a function of $\kappa$ it is locally bounded.
Thus, the sequence
$\{v_k\}$ contains
a subsequence, again denoted by $\{v_k\}$, that converges locally uniformly
to a bounded Lipschitzian function
$v\in W_{loc}^{2,2}(\mathbb{R}^d)$ such that
$$
\sup\limits_x|v(x)|
\leq \lambda^{-1} \sup\limits_x|f(x)|
\quad\hbox{and}\quad
\sup\limits_x|\nabla v(x)|
\leq \lambda^{-1} \sup\limits_x|\nabla f(x)|,
$$
and, in addition,  the restrictions
of $v_k$ to any ball $U$
converge  to~$v|_U$ weakly in $W^{2,2}(U)$.

Let $\widehat{L}$ be the elliptic operator with
the same second order part as $L$,
but with
drift is $\widehat{b}=2A\nabla \varrho/\varrho -b$.
Then by the integration by parts formula
$$
\int \psi L\varphi \, d\mu=
\int \varphi \widehat{L}\psi\, d\mu
\quad \hbox{for all $\psi ,\varphi\in C_0^\infty(\mathbb{R}^d)$.}
$$
In addition, for any $\lambda>0$,
the ranges of $\lambda -L$ and
$\lambda - \widehat{L}$ on $C_0^\infty(\mathbb{R}^d)$ are dense
in $L^1(\mu)$.
The operator $\widehat{L}$ also generates a Markov semigroup
on $L^1(\mu)$ with respect to which $\mu$ is
invariant. The corresponding resolvent is denoted
by~$\widehat{G}_\lambda$.
For the proofs we refer to
\cite[Proposition 2.9]{BRS0} or
\cite[Proposition 1.10(b)]{Stannat}
(see also \cite[Theorem 3.1]{BRS}).

Now we show that $v=G_\lambda f$.
Note that $\varrho_k\to\varrho$ uniformly on balls
according to \cite{BR},~\cite{BKR}.
Hence, given $\varphi\in C_0^\infty(\mathbb{R}^d)$ with support
in a ball~$U$,
we have
$$
\int [\lambda v - Lv -f]\varphi\varrho\, dx
=\lim\limits_{k\to\infty}
\int [\lambda v_k - L_kv_k -f]\varphi\varrho_k\, dx=0
$$
by weak convergence of $v_k$ to $v$ in $W^{2,2}(U)$
combined with convergence of $b_k$ to $b$ in
$L^p(U,\mathbb{R}^d)$.
Therefore, by the integration by parts formula
$$
\int v(\lambda\varphi -\widehat{L}\varphi)\, d\mu=\int f\varphi\,
d\mu
$$
for all $\varphi\in C_0^\infty(\mathbb{R}^d)$. The function
$G_\lambda f$ is bounded and
satisfies the same relation, so it remains to recall
that if a bounded function $u$ satisfies
$\int u(\lambda\varphi -\widehat{L}\varphi)\, d\mu=0$ for
all $\varphi\in C_0^\infty(\mathbb{R}^d)$, then $u=0$~a.e.,
since $(\lambda -\widehat{L})\bigl(C_0^\infty(\mathbb{R}^d)\bigr)$
is dense in~$L^1(\mu)$.

Now we turn to the pointwise estimate
     $|\nabla G_\lambda f(x)|\le \lambda^{-1}G_\lambda|\nabla f|(x)$.
Suppose first that $f\in C_0^\infty(\mathbb{R}^d)$.
The desired estimate holds for every $G_\lambda^{(k)}$ in place
     of~$G_\lambda$. It has been shown above that
     $v=G_\lambda f$ is a weak limit of $v_k=G_\lambda^{(k)}f$
     in $W^{2,2}(U)$ for every ball~$U$.
In addition, the functions $G_\lambda^{(k)}|\nabla f|$ converge
     weakly in $W^{2,2}(U)$ to the function $G_\lambda|\nabla f|$, which is
     also clear by the above reasoning.
     Since the embedding of $W^{2,2}(U)$ into $W^{2,1}(U)$ is compact,
     we may assume, passing to a subsequence, that
    $\nabla G_\lambda^{(k)}f(x)\to \nabla G_\lambda f(x)$
and $G_\lambda^{(k)}|\nabla f|(x)\to G_\lambda |\nabla f|(x)$
almost everywhere on~$U$. Hence we arrive at the desired estimate.
If $f$ is Lipschitzian and has bounded support,
we can find uniformly Lipschitzian functions
     $f_n\in C_0^\infty(\mathbb{R}^d)$ vanishing outside some
     ball such that $f_n\to f$ uniformly and
     $\nabla f_n\to \nabla f$~a.e. Then, by the same reasons as above,
     one has
     $G_\lambda |\nabla f_n|\to G_\lambda |\nabla f|$ and
    $\nabla G_\lambda f_n\to \nabla G_\lambda f$ in $L^2(U)$.
     Passing to an almost everywhere convergent subsequence
     we obtain a pointwise inequality. Finally,
     in the case of a general Lipschitzian function $f\in L^1(\mu)$,
we can find uniformly Lipschitzian functions
     $\zeta_n$ such that $0\le\zeta_n\le 1$ and $\zeta_n(x)=1$ if
     $|x|\le n$. Let $f_n=f\zeta_n$. By the previous step
     we have
$$
|\nabla G_\lambda f_n(x)|\le \lambda^{-1}G_\lambda|\nabla f_n|(x).
$$
     The functions $f_n$ are uniformly Lipschitzian. Hence, for every
     ball~$U$, the sequence of functions $G_\lambda f_n|_U$ is
          bounded in the norm of
     $W^{2,2}(U)$. In addition, the functions
$G_\lambda |\nabla f_n|$ on $U$ converge to $G_\lambda |\nabla f|$ in
$L^2(U)$, since $|\nabla f_n|\to |\nabla f|$ in $L^2(\mu)$
by the Lebesgue dominated convergence theorem.
     Therefore, the same reasoning as above
     completes the proof.
\end{proof}

\begin{proof}[Proof of Theorem~2]
Suppose first that $A$ is piece-wise constant, i.e.,
there exist finitely many intervals $[0,t_1)$, $[t_1,t_2)$,\ldots,
$[t_{n},1]$ such that $A(t)=A_k$ whenever $t_{k-1}\le t<t_k$,
where each $A_k$ is a strictly positive symmetric matrix.
In addition, let us assume that there exist vector fields $b_k$
such that $b(t,x)=b_k(x)$ whenever $t_{k-1}\le t<t_k$.
Then we obtain a solution
$u$ by successively applying the semigroups
$T_t^{(k)}$ generated by the elliptic operators
with the diffusion matrices $A_k$ and drifts $b_k$, i.e.,
$$
u(t,x)=T_{t-t_{k-1}}T_{t_{k-1}}\cdots T_{t_1}f(x)\quad
\hbox{whenever
$t\in [t_{k-1},t_k)$.}
$$
The conclusion of Theorem~2 in this case
follows by Theorem~1. Our next step is to approximate $A$ and $b$
by mappings of the above form in such a way that the corresponding
sequence of solutions would converge to a solution of our equation.
Let us observe that, for an arbitrary sequence of such solutions
$u_k$ corresponding to piece-wise constant in time coefficients,
for every compactly supported function $\varphi$ on $\mathbb{R}^d$,
the functions
\begin{equation}\label{ulip}
t\mapsto \int_{\mathbb{R}^d} \varphi(x)u_k(t,x)\, dx
\end{equation}
are uniformly Lipschitzian provided that the operator norms
of the matrix functions $A_k$ are uniformly bounded and that
the $L^p(K)$-norms of the vector fields $b_k(t,\,\cdot\,)$
are uniformly bounded for every fixed cube $K$ in $\mathbb{R}^d$.
This is clear, because (\ref{par1}) can be written as
$$
\int_{\mathbb{R}^d}\varphi(x)u(t,x)\, dx =
\int_0^t\int_{\mathbb{R}^d}
[L_s\varphi(x) \, u(s,x)
+\varphi(x)b^i(s,x)\partial_{x_i}u(s,x)]\, dx \, ds ,
$$
where in the case $u=u_k$ we have
$$
|u(s,x)|\le \sup |f(x)|
\quad\hbox{and}\quad
|\nabla_x u(s,x)|\le \sup |\nabla f(x)|.
$$
One can choose a subsequence in $\{u_k\}$ that converges
to some function $u$ on $[0,1]\times\mathbb{R}^d$ in the following
sense: for every cube $K$ in $\mathbb{R}^d$, the functions
the restrictions of the functions $u_k$ to $[0,1]\times K$
converge weakly to $u$ in the space $L^2([0,1],W^{2,2}(K))$,
where each $u_k$ is regarded as a mapping
$t\mapsto u_k(t,\,\cdot\,)$ from $[0,1]$ to $W^{2,2}(K)$.
Passing to another subsequence we obtain
$$
\lim\limits_{n\to\infty}
\int_{\mathbb{R}^d}\varphi(x)u_k(t,x)\, dx=
\int_{\mathbb{R}^d}\varphi(x)u(t,x)\, dx
$$
for all $t\in [0,1]$ and
all smooth compactly supported~$\varphi$.
Indeed, for a given function $\varphi$ this is possible due to the uniform
Lipschitzness of the functions~(\ref{ulip}).
Then our claim is true for a countable
family of functions~$\varphi$, which, on account of the uniform
boundedness of $u_k$, yields the claim for all~$\varphi$.
Therefore, it remains to find approximations $A_k$ and $b_k$
such that, for every function $\psi\in C[0,1]$, the integrals
$$
\int_0^1\psi(s)\int_{\mathbb{R}^d}[L_s^{(k)}\varphi(x) \, u_k(s,x)
+\varphi(x)b^i_k(s,x)\partial_{x_i}u_k(s,x)]\, dx \, ds
$$
would converge to the corresponding integral with $A,b$, and~$u$.
Clearly, it suffices to obtain the desired convergence
for suitable countable families of functions
$\varphi_i$ and~$\psi_j$. Let us fix two sequences
$\{\psi_j\}\subset C[0,1]$ and $\{\varphi_i\}\subset
C_0^\infty (\mathbb{R}^d)$ with the following property:
every compactly supported square-integrable function
$v$ on $[0,1]\times\mathbb{R}^d$ can be approximated
in $L^2$ by a sequence of finite linear combinations
of products $\psi_j\varphi_i$.
Let us consider the functions
$$
\alpha_{i,j,k}(t):=a^{ij}(t)\psi_k(t),\quad
\beta_{i,j,k}(t):=
\psi_k(t)\int_{\mathbb{R}^d} b^i(s,x)\varphi_j(x)\, dx,
$$
$$
\theta_{k,i}(t)=\int_{[-k,k]^d} b_i(t,x)^2\, dx.
$$
Let $\mathcal{F}$
denote the obtained countable family of functions extended
periodically from $[0,1)$ to $\mathbb{R}$ with period~$1$.
It is well known that, for almost every $s\in [0,1)$,
the Riemannian sums $R_n(\theta)(s)=2^{-n}\sum\limits_{k=1}^{2^n}
\theta(s+k2^{-n})$ converge to the integral of $\theta$ over
$[0,1]$ for each $\theta\in\mathcal{F}$.
It follows that one can find points $t_{n,l}$, $l=1,\ldots,N_n$,
$n\in\mathbb{N}$, such that
$$
0=t_{n,0}<t_{n,1}<t_{n,2}<\cdots<t_{n,N_n}=1
$$
and, for every
$\theta\in \mathcal{F}$, letting
$\theta_n(t):=\theta(t_{n,l})$ whenever
$t_{n,l-1}\le t<t_{n,l}$, one has
$$
\int_0^1 \theta_n(t)\, dt\to\int_0^1 \theta(t)\, dt .
$$
To this end, we pick a common point $s_0$ of convergence of the
Riemann sums $R_n(\theta)(s_0)$ to the respective integrals and
let $t_{n,l}=s_0+l2^{-n}$ (mod$1$).
By using the points $t_{n,l}$, one obtains the desired
piece-wise constant approximations of $A$ and~$b$.
Namely, let
$A_n(t)=A(t_{n,l})$ and $b_n(t,x)=b(t_{n,l},x)$ whenever
$t_{n,l-1}\le t< t_{n,l}$.
As explained above,
passing to a subsequence, we may assume that
the corresponding solutions $u_n$ converge to a function
$u$ such that, for every cube $K=[-m,m]^d$ in~$\mathbb{R}^d$ and every
$t\in (0,1]$, one has
$$
u(t,\,\cdot\,)|_K\in W^{2,2}(K),\quad
\int_0^1 \|u(t,\,\cdot\,)\|_{W^{2,2}(K)}^2\, dt<\infty ,
     $$
and for any function $\zeta\in L^2([0,1]\times K)$ there holds the
equalities
\begin{align*}%%\label{twoeq}
\lim\limits_{n\to\infty}
\int_0^1\int_K \zeta(t,x)u_n(t,x)\, dx\, dt
&=
\int_0^1\int_K \zeta(t,x)u(t,x)\, dx\, dt ,
\\
\lim\limits_{n\to\infty}
\int_0^1\int_K \zeta(t,x)\partial_{x_i}\partial_{x_j}u_n(t,x)\, dx\, dt
&=
\int_0^1\int_K \zeta(t,x)
\partial_{x_i}\partial_{x_j}u(t,x)\, dx\, dt ,
\\
\lim\limits_{n\to\infty}
\int_0^1\int_K
\zeta(t,x)\partial_{x_i}u_n(t,x)\, dx\, dt
&=
\int_0^1\int_K \zeta(t,x)
\partial_{x_i}u(t,x)\, dx\, dt ,
\\
\lim\limits_{n\to\infty}
\int_0^1\int_K
b_i^n(t,x)^2\, dx\, dt
&=
\int_0^1\int_K
b_i(t,x)^2\, dx\, dt .
\end{align*}
Note that for every cube $K$ in $\mathbb{R}^d$, the restrictions
of the functions $b_n^i$ to $[0,1]\times K$ converge
to the restriction of $b^i$ in the norm of $L^2([0,1]\times K)$.
This is clear from the last displayed equality, which gives
convergence of $L^2$-norms, along with convergence
of the Riemann sums $R_n(\beta_{i,j,k})(s_0)$ to the integral of
$\beta_{i,j,k}$ over $[0,1]$, which yields weak convergence
(we recall that if a sequence of vectors $h_n$ in a Hilbert space
$H$ converges weakly to a vector $h$ and the norms of $h_n$ converge
to the norm of~$h$, then there is norm convergence).
It follows that for any $\psi\in C[0,1]$ and any
$\varphi\in C_0^\infty(\mathbb{R}^d)$ with support in
$[-m,m]^d$, we have
\begin{multline*}
\lim\limits_{n\to\infty}
\int_0^1
\psi(t)a_n^{ij}(t)
\int_{\mathbb{R}^d}
\partial_{x_i}\partial_{x_j}\varphi(x)u_n(t,x)\, dx\, dt
\\
=\int_0^1
\psi(t)a^{ij}(t)
\int_{\mathbb{R}^d}
\partial_{x_i}\partial_{x_j}\varphi(x)u(t,x)\, dx\, dt.
\end{multline*}
In addition,
\begin{multline*}
\lim\limits_{n\to\infty}
\int_0^1
\psi(t)\int_{\mathbb{R}^d}
\varphi(x)\partial_{x_i}u_n(t,x)b_n^i(t,x)\, dx\, dt
\\
=\int_0^1
\psi(t)
\int_{\mathbb{R}^d}
\varphi(x)\partial_{x_i}u(t,x)b^i(t,x)\, dx\, dt .
\end{multline*}
This follows by norm convergence of
$b_n^i$ to $b^i$
and
weak convergence of
$\varphi\partial_{x_i}u_n$ to
$\varphi\partial_{x_i}u$
in~$L^2([0,1]\times [-m,m]^d)$.
Therefore,
for every $\varphi\in C_0^\infty(\mathbb{R}^d)$,
one has
$$
\int_{\mathbb{R}^d}
\varphi(x)u(t,x)\, dx\, dt
=
\int_{\mathbb{R}^d}
\varphi(x)f(x)\, dx+
\int_0^t\int_{\mathbb{R}^d}
\varphi(x)L_tu(t,x)\, dx\, dt
$$
for almost all $t\in [0,1]$, since the integrals of both sides
multiplied by any function $\psi\in C_0^\infty(0,1)$ coincide.
Taking into account the continuity of both sides (the left-hand
side is Lipschitzian as explained above), we conclude that
the equality holds for all $t\in [0,1]$.
\end{proof}

{\bf Acknowledgements.}
{\small
This work has been supported in part by
the RFBR project
04-01-00748,
the DFG Grant 436 RUS 113/343/0(R),
the INTAS project 03-51-5018,
the Scientific Schools
Grant 1758.2003.1,
the DFG--Forschergruppe ``Spectral Analysis, Asymptotic
Distributions, and Stochastic Dynamics'',
the BiBoS--reseach centre,
and the research programme ``Analisi e controllo
di equazioni di evoluzione deterministiche e
stocastiche'' from the Italian ``Ministero della
Ricerca Scientifica e Tecnologica''.}

%%V.B.: Department of Mechanics and Mathematics,
%%  Moscow State University, 119992 Moscow, Russia,

%%G.D.: Scuola Normale Superiore di Pisa, Piazza dei Cavalieri 7,
%%I--56125 Pisa, Italy,

%%M.R.: Fakult\"at f\"ur Mathematik,
%%Universit\"at Bielefeld, D--33501 Bielefeld, Germany,

%%Z.S.: Department of Mathematics, Imperial College,
%%180 Queens Gate, London SW7 2BZ, UK

\end{document}